\documentclass[a4paper]{amsart}\usepackage{amsfonts,amsmath,amssymb,bbm}

\def\k{\mathbbm{k}}

\def\R{\Bbb{R}}\def\Z{\Bbb{Z}}

\def\li{\ \\ $\bullet$ }

\def\suml{\sum\limits}\def\coprodl{\coprod\limits}

\def\prodl{\mathop\prod\limits}

\newcommand{\quotients}[2]{{\footnotesize\left.\raisebox{0.4ex}{$#1$}\! / \!\raisebox{-0.4ex}{$#2$}\right.}}

\def\tg{\tilde{g}}

\def\al{\alpha}\def\ga{\gamma}\def\Ga{\Gamma}\def\be{\beta}

\def\la{\lambda}\def\si{\sigma}

\def\cK{{\mathcal K}}

\def\uA{{\underline{A}}}
\def\u0{\underline{0}}

\def\uk{{\underline{k}}}

\def\ux{{\underline{x}}}

\def\one{{1\hspace{-0.1cm}\rm I}}

\newcommand{\ber}{\begin{array}{l}}\newcommand{\eer}{\end{array}}
\newcommand{\bpm}{\begin{pmatrix}}\newcommand{\epm}{\end{pmatrix}}
\newcommand{\bM}{\begin{matrix}}\newcommand{\eM}{\end{matrix}}
\newcommand{\bee}{\begin{enumerate}}\newcommand{\eee}{\end{enumerate}}

\def\sset{\subset}\def\sseteq{\subseteq}\def\smin{\setminus}

\def\ND{Newton diagram}

\def\bull{\vrule height .9ex width .9ex depth -.1ex }

\newcommand{\beq}{\begin{equation}}\newcommand{\eeq}{\end{equation}}
\newcommand{\beqm}{\begin{multline}}\newcommand{\eeqm}{\end{multline}}
\newtheorem{Lemma}{Lemma}[section]\newcommand{\bel}{\begin{Lemma}}\newcommand{\eel}{\end{Lemma}}
\newtheorem{Example}[Lemma]{Example}\newcommand{\bex}{\begin{Example}\rm}\newcommand{\eex}{\end{Example}}
\newtheorem{Proposition}[Lemma]{Proposition}\newcommand{\bprop}{\begin{Proposition}}\newcommand{\eprop}{\end{Proposition}}
\newtheorem{Definition-Proposition}[Lemma]{Definition-Proposition}

\def\bpr{~\\{\em Proof.\ }}\newcommand{\epr}{$\bull$\\}
\newtheorem{Theorem}[Lemma]{Theorem}\newcommand{\bthe}{\begin{Theorem}}\newcommand{\ethe}{\end{Theorem}}
\newtheorem{Definition}[Lemma]{Definition}\newcommand{\bed}{\begin{Definition}}\newcommand{\eed}{\end{Definition}}
\newtheorem{Remark}[Lemma]{Remark}\newcommand{\beR}{\begin{Remark}\rm}\newcommand{\eeR}{\end{Remark}}
\newtheorem{Corollary}[Lemma]{Corollary}\newcommand{\bcor}{\begin{Corollary}}\newcommand{\ecor}{\end{Corollary}}
\newcommand{\bet}{\begin{tabular}{cccccccc}}\newcommand{\eet}{\end{tabular}}

\newcommand{\bin}[2]{\binom{#1}{#2}}
\newcommand\isom{\xrightarrow{\,\smash{\raisebox{-0.65ex}{\ensuremath{\scriptstyle\sim}}}\,}}

\title[]{A \MakeLowercase{generalized} FKG-\MakeLowercase{inequality for compositions}}
\author[]{D\MakeLowercase{mitry} K\MakeLowercase{erner} \MakeLowercase{and} A\MakeLowercase{ndr\'as} N\MakeLowercase{\'emethi}}
\address{
Department of Mathematics, Ben Gurion University of the Negev, P.O.B. 653, Be'er Sheva 84105, Israel,
dmitry.kerner@gmail.com,\newline
R\'enyi Institute of Mathematics, Budapest, Re\'altanoda u. 13--15, 1053, Hungary,
nemethi@renyi.hu.}
\date{\today}
\thanks{D.K. was supported by the grant FP7-People-MCA-CIG, 334347.
A.N. was supported by  OTKA Grant 100796.
\\
We thank Anders Bj\"{o}rner, Efim Dinitz, P\'{a}lv\"{o}lgyi D\"{o}m\"{o}t\"{o}r, L\'{a}szl\'{o}
 Lov\'{a}sz, Ron Peled,  Siddhartha Sahi, Michael Saks, Gjergji Zaimi for valuable advices.}

\begin{document}
\begin{abstract}
We prove a Fortuin-Kasteleyn-Ginibre-type inequality for the lattice of compositions of
 the integer $n$ with at most $r$ parts. As an immediate application we get a wide  generalization of the classical Alexandrov-Fenchel inequality for mixed volumes and of Teissier's inequality for mixed covolumes.
\end{abstract}
\maketitle \setcounter{secnumdepth}{6} \setcounter{tocdepth}{1}

\section{Introduction}
\subsection{} Consider a finite partially ordered set $(X,\preceq)$ and two non-decreasing (non-negative) functions, $f,g:X\to\R_{\ge0}$.
(Namely, for any $x,y\in X$, if $x\preceq y$ then one has $f(x)\le f(y)$ and $g(x)\le g(y)$.)
The product function $f\cdot g:\ X\to\R_{\ge0}$ is also non-decreasing.
Take the arithmetic average $$Av_X(f):=(\textstyle{\suml_{x\in X}}f(x))/|X|.$$
  A natural question is whether  $Av_X(f)\cdot Av_X(g)$ can be compared with $Av_X(f\cdot g)$.
\bex
Suppose that $X$ is totally ordered. Then the non-decreasing functions are just the non-decreasing sequences
 of real numbers, $0\le a_1\le\cdots\le a_n$ and $0\le b_1\le\cdots\le b_n$.
 In this case
 the comparison of the averages is realized by
  the classical Chebyshev sum inequality: $(\sum_i a_i)(\sum_j b_j)\le n(\sum_i a_ib_i)$.
\eex
On the other hand, if the order on $X$ is not ``strong enough'' then the inequality utterly fails. Hence, the more precise question is:
\beq\label{Eq.Main.Question}
\text{Which posets  does }Av_X(f)\cdot Av_X(g)\le Av_X(f\cdot g) \text{ hold for}?
\eeq
If $(X,\preceq)$ admits an action of some group $G$, then one can consider the ``equivariant'' version of this question
by taking  $G$-invariant functions $f$ and $g$.

The fundamental Fortuin-Kasteleyn-Ginibre (FKG) inequality settles the question for a large class of lattices:
\bthe\cite{FKG-1971}, see also \cite[pg. 147, Theorem 5]{Bollobas}.
\\Let $X$ be a finite distributive lattice. Consider a ``measure'', $X\stackrel{\mu}{\to}\R_{\ge0}$, which is log-supermodular, i.e. $\mu(x\wedge y)\mu(x\vee y)\ge\mu(x)\mu(y)$ for any $x,y\in X$.
  Then $\Big(\suml_{x\in X}f(x)g(x)\mu(x)\Big)\cdot \suml_{x\in X}\mu(x)\ge\Big(\suml_{x\in X}f(x)\mu(x)\Big)\Big(\suml_{x\in X}g(x)\mu(x)\Big)$.
\ethe
(The inequality of equation \eqref{Eq.Main.Question} is obtained for the constant measure, $\mu(x)=1$, which is trivially supermodular.)

One of the interpretation of the FKG inequality is: ``in many systems the increasing events are positively correlated'' (while an increasing event and a decreasing event are negatively correlated).
Hence, the applications of this inequality go far beyond the combinatorics and include e.g.  statistical mechanics and probability.

\

\subsection{}
The condition ``$X$ is a distributive lattice'' in the above theorem
is rather restrictive. Many of the natural posets appearing in arithmetics/algebra/geometry are not of this type.
In the current work we establish the inequality of equation  \eqref{Eq.Main.Question} for a particular poset $\cK_{n,r}$ of ordered compositions, cf. Theorem \ref{Thm.Averaging.over.Knr}. This poset appears frequently  in the context of the Young diagrams
(representation theory), complete intersections (algebraic geometry), mixed (co)volumes/multiplicities
  (integral geometry and commutative algebra).

It is known that lattices for which the FKG inequality holds (for any log-supermodular measure) are necessarily distributive. Theorem \ref{Thm.Averaging.over.Knr} gives an example of non-distributive lattice for which FKG holds for a particular measure.

For related inequalities see also \cite{Cu.Gr.Sk}.

\subsection{} Sometimes one considers the geometric average, $Av^G_X(f):=(\prod_{x\in X}f(x))^{\frac{1}{|X|}}$
as well. Then it is natural to compare $(Av^G_X(f))^{Av_X(g)}$
 to $Av^G_X(f^g)$. For example, given the real numbers $1\le a_1\le \cdots\le a_n$ and $0\le b_1\le\cdots\le b_n$, one has: $(\sqrt[n]{\prod a_i})^{\frac{\sum b_i}{n}}\le\sqrt[n]{\prod a_i^{b_i}}$.

One passes between $Av_X$ and $Av^G_X$ by using $ln$ and $exp$. Thus a simple reformulation of Theorem \ref{Thm.Averaging.over.Knr} provides automatically  the comparison of the geometric averages as well,
 cf. Corollary \ref{Thm.Inequality.Geometric.Averages}.

\subsection{}
As an immediate application in \S\ref{Sec.Alexandrov.Fenchel.Teissier.ineqs}
we prove  a highly non-trivial convexity property for the mixed volumes of convex bodies in $\R^N$
 and the parallel statement for the mixed volumes of \ND s (i.e. co-volumes of the Newton polyhedra).
 These generalize the classical Alexandrov-Fenchel inequality and its Teissier's counterpart. In particular, it gives a partial answer to a question of \cite{gromov1990convex}.

Furthermore, in \cite{Kerner-Nemethi.Durfee.Nnd} we heavily use the inequality of averages to establish a bound on some topological invariants of singularities.
In fact this was our initial motivation for the inequality of averages.

\section{The poset $\cK_{n,r}$}
\subsection{The set of compositions}\label{Sec.Knr.first.time}
Denote by $\cK_{n,r}$ the set of the (ordered) compositions of the integer $n$ into $r$ summands,
\beq
\cK_{n,r}:=\{\uk=(k_1,\ldots, k_r)\, :\,  k_i\geq 0 \ \mbox{for all $i$, and } \ \sum _ik_i=n\}.
\eeq
This $\cK_{n,r}$ can be thought of as the lattice points $\uk$ of the simplex
$\{\sum_ix_i=n\}\cap {\mathbb R}^r_{\geq 0}$.
 Its  cardinality is $|\cK_{n,r}|=\binom{n+r-1}{n}$.

The permutation group of $r$ elements, $\Xi_r$, acts on $\cK_{n,r}$ by $\si(\uk)=(k_{\si(1)},\dots,k_{\si(r)})$.  The quotient $\quotients{\cK_{n,r}}{\Xi_r}$ is the set of partitions into $r$ summands.
(In other words,  a partition is an {\em unordered} composition.)
For the general introduction see \cite[Chapter 7]{Stanley}.

For convenience we put $\cK_{n,r}=\varnothing$ when $r\le0$ or $n<0$.

\subsection{}{\bf The `dominance'  order on $\cK_{n,r}$} can be defined as follows.
\li Suppose for $\uk=(k_1,k_2,\dots,k_r)\in\cK_{n,r}$ one has $k_1-1\ge k_2+1$.
Then put $\uk\succeq \uk':=(k_1-1,k_2+1,k_3,\dots,k_r)$.
\li Extend this by transitivity, i.e. if $\uk\succeq\uk'$ and $\uk'\succeq\uk''$ then $\uk\succeq\uk''$.
\li Extend this by the action $\Xi_r\circlearrowright\cK_{n,r}$, i.e. if $\uk\succeq\uk'$ then $\si(\uk)\succeq\si(\uk')$ for
any $\si\in\Xi_r$.

In this way we get a partially ordered set with the maximal elements, $(n,0,\dots,0)$ and its orbit under $\Xi_r$,
and the minimal elements,
$(\lfloor\frac{n}{r}\rfloor,\dots,\lfloor\frac{n}{r}\rfloor,\lceil\frac{n}{r}\rceil,\dots,\lceil\frac{n}{r}\rceil)$
 and its orbit under $\Xi_r$.
 By construction, this partial order is $\Xi_r$ invariant.
 In particular, any two different elements inside a $\Xi_r$ orbit are incomparable.
  This order descends to the quotient $\quotients{\cK_{n,r}}{\Xi_r}$. (Indeed, for any two elements of $\quotients{\cK_{n,r}}{\Xi_r}$,
  suppose some of their representatives in $\cK_{n,r}$ are comparable, $\uk\preceq\uk'$.
  Then put $[\uk]\preceq [\uk']$. By the $\Xi_r$-invariance this assignment is consistent: if some other
   preimages are comparable then they satisfy the same inequality.)
    The quotient poset $\quotients{\cK_{n,r}}{\Xi_r}$ has a unique minimal and a unique maximal element.
 \bex
1. For $\cK_{n,2}$ we get $(n,0)\succeq(n-1,1)\succeq\cdots\succeq(\lceil\frac{n}{2}\rceil,\lfloor\frac{n}{2}\rfloor)$
 and $(\lfloor\frac{n}{2}\rfloor,\lceil\frac{n}{2}\rceil)\preceq\cdots\preceq(0,n)$.
 For $\quotients{\cK_{n,2}}{\Xi_2}$: $(n,0)\succeq(n-1,1)\succeq\cdots\succeq(\lceil\frac{n}{2}\rceil,\lfloor\frac{n}{2}\rfloor)$.
In particular, $\quotients{\cK_{n,2}}{\Xi_2}$ is {\em totally} ordered.
\\2. As mentioned above, $\cK_{n,r}$ is never totally ordered for $r>1$:
 the elements of any $\Xi_r$-orbit are incomparable.
The quotient $\quotients{\cK_{n,r}}{\Xi_r}$ is not totally ordered for $r\ge3$ and high enough $n$.
 For example, for $(n,r)=(6,3)$, the elements $(4,1,1)$ and $(3,3,0)$ are incomparable.
\\3. In the particular case $\quotients{\cK_{n,n}}{\Xi_n}$ coincides with the Young lattice of all the possible partitions of $n$, 
see e.g. \cite{Brylawski}, \cite[Chapter 7]{Stanley}.
\\4. $\quotients{\cK_{n,r}}{\Xi_r}$ is always a lattice, but it is non-distributive
for $n\ge7$. It contains the non-distributive sublattice:
\beq
\bM \overset{(*,4,2,1)}{\bullet}&\to&\overset{(*,4,1,1,1)}{\bullet}&
\\\downarrow&&&\searrow
\\\underset{(*,3,3,1)}{\bullet}&\to&\underset{(*,3,2,2)}{\bullet}&\to&\underset{(*,3,2,1,1)}{\bullet}\eM.
\eeq
Here $*$ denotes the rest of the partition, some fixed tuple whose sum is $(n-7)$ and such that 
 all the entries go in the non-increasing order. 
For this lattice one has:
\begin{multline}
(*,3,3,1)\wedge\Big((*,4,1,1,1)\vee(*,3,2,2)\Big)=(*,3,3,1),
\\
\Big((*,3,3,1)\wedge(*,4,1,1,1)\Big)\vee\Big((*,3,3,1)\wedge(*,3,2,2)\Big)=(*,3,2,2).
\end{multline}
\eex

\subsection{}
Suppose that a set of objects is indexed by this set of compositions, $\{A_\uk\}_{\uk\in\cK_{n,r}}$.
 We often use the standard set-theoretical inclusion-exclusion formula:
\beq\label{Eq.Exclusion.Inclusion.Formula}
\suml_{\uk\in\cK_{n,r}}A_\uk-\suml^r_{i=1}\suml_{\substack{\uk\in\cK_{n,r}\\k_i=0}}A_\uk+
\suml_{1\le i_1<i_2\le r}\suml_{\substack{\uk\in\cK_{n,r}\\k_{i_1}=k_{i_2}=0}}A_\uk-\cdots=
\suml_{\substack{\uk\in\cK_{n,r}\\k_1,\dots,k_r>0}}A_\uk.
\eeq

\section{The inequality for  averages over $\cK_{n,r}$}\label{Sec.Convexity.Combinatorial}
\bthe\label{Thm.Averaging.over.Knr}
Let $f,g:\cK_{n,r}\to\R_{\ge0}$ be non-negative functions. Suppose $f$ is $\Xi_r$-invariant,
 i.e. $f(\si(x))=f(x)$ for any $\si\in\Xi_r$.
\\(a) If both functions are non-decreasing then $Av_{\cK_{nr}}(f)\cdot Av_{\cK_{n,r}}(g)\le Av_{\cK_{n,r}}(fg)$.
\\(b) If $f$ is non-decreasing, while $g$ is non-increasing
 then $Av_{\cK_{n,r}}(f)\cdot Av_{\cK_{n,r}}(g)\ge Av_{\cK_{n,r}}(fg)$.
\\(c) In the statements  above, the equality holds if and only if  either $f$ is constant, or the symmetrization of $g$,
$\tg(x):=Av_{\Xi_r(x)}(g)$ is constant on $\cK_{n,r}$.
\ethe
\bex
1. This inequality was proved in \cite[\S1.3]{Kerner-Nemethi.Durfee.homogen}
in the particular case of
 $f(\uk):=\prodl^r_{i=1}d^{k_i}_i$ and $g(\uk):=\bin{n+r}{k_1+1,\dots,k_r+1}$, in which
  case it reads as
\beq\label{Eq.Combinatorial.Inequality}
\Big(\suml_{\uk\in\cK_{n,r}}\bin{n+r}{k_1+1,\dots,k_r+1}\Big)\Big(\suml_{\uk\in\cK_{n,r}}(\prod^r_{j=1}d^{k_j}_j)\Big)
\ge|\cK_{n,r}|\suml_{\uk\in\cK_{n,r}}\bin{n+r}{k_1+1,\dots,k_r+1}\Big(\prod^r_{j=1}d^{k_j}_j\Big).
\eeq
This was the essential ingredient in establishing a bound for some topological invariants in Singularity Theory.

2. If one does not assume that at least one of the functions is $\Xi_r$-invariant then the inequality does not hold.
For example, consider $\cK_{2,2}=\{(0,2) \succeq (1,1) \preceq (2,0)\}$. The inequality is obviously violated for
\[f(x)=\Big\{\ber 1,\ x=(0,2)\\0,\ otherwise\eer,\quad\quad g(x)=\Big\{\ber 1,\ x=(2,0)\\0,\ otherwise\eer.\]

3. The statement can be reformulated as averaging over $\quotients{\cK_{n,r}}{\Xi_r}$, with the weight function $\mu([\uk])=|\Xi_r(\uk)|$.
\eex
\bpr
{\bf Step 1.} We start with some simplifying remarks.

We can (and will) assume that  $g$ is  $\Xi_r$-invariant. Indeed, define $\tg:\cK_{n,r}\to\R_{\ge0}$
by the averaging over the orbit, $\tg(x):=Av_{\Xi_r(x)}(g)$.
By construction $\tg$ is non-negative, non-decreasing and $\Xi_r$-invariant.
 And $Av_{\cK_{n,r}}(g)=Av_{\cK_{n,r}}(\tg)$, $Av_{\cK_{n,r}}(fg)=Av_{\cK_{n,r}}(f\tg)$.

\

Note that the statement (b) follows from (a).
 (Indeed, if $g$ is non-increasing then consider the number $\underset{\cK_{n,r}}{max}(g)$ and apply the first
 statement to the functions $f$ and $\underset{\cK_{n,r}}{max}(g)-g$.)
In particular, we can assume that both $f$ and $g$ are non-decreasing and $\Xi_r$-invariant functions.

Suppose $(f_1,g)$ and $(f_2,g)$ satisfy the statement of the
theorem. Then the pair $(a_1 f_1+a_2f_2,g)$ satisfies the theorem  too
 for any $a_1,a_2\in\R_{\ge0}$. (Note that $a_1 f_1+a_2f_2$ is $\Xi_r$-invariant and non-decreasing as well.)
 Moreover, the equality holds iff it holds for each pair separately: $Av_{\cK_{n,r}}(f_i)Av_{\cK_{n,r}}(g)=Av_{\cK_{n,r}}(f_i g)$. On the other hand,
  any $\Xi_r$-invariant, non-decreasing function on $\cK_{n,r}$ is presentable as a positive linear combination
  of certain  $\{0,1\}$ valued functions, which  are $\Xi_r$-invariant and non-decreasing as well.
  Therefore, to prove (a),  we can assume that both $f$ and $g$ are of this form.

\

Therefore, we assume that $f$ and $g$ are the characteristic functions of some subsets. E.g.,  $f=\one_X$, where $X\sseteq\cK_{n,r}$ is $\Xi_r$-invariant
and ``upward closed''. This later condition  means the following:
 if $\uk\in X$ and $\uk'\succeq\uk$ then $\uk'\in X$.
First and last,  it is enough to prove the inequality for
cardinalities of ``upward closed''
and $\Xi_r$-invariant
subsets:
\beq
Av_{\cK_{n,r}}(X)\cdot Av_{\cK_{n,r}}(Y):=\frac{|X|}{|\cK_{n,r}|}\cdot \frac{|Y|}{|\cK_{n,r}|}\le \frac{|X\cap Y|}{|\cK_{n,r}|}=Av_{\cK_{n,r}}(X\cap Y).
\eeq
Furthermore, we also verify that the equality holds iff either $X=\cK_{n,r}$ or $X=\varnothing$.

 Our proof generalizes and refines
 the  proof of \cite[\S4]{Kerner-Nemethi.Durfee.homogen}.

\

{\bf Step 2.} The proof below consists of some counting over the subsets of $\cK_{n,r}$. First we define the stratification:
 $\cK_{n,r}=\coprodl_{s=0,\dots,r-1}\cK^s_{n,r}$, where $\cK^s_{n,r}:=\{\uk\in\cK_{n,r}:\ |\{i:\ k_i=0\}|=s\}$.
 Note that $\cK^s_{n,r}=\varnothing$ for $s<r-n$.
 We often use the expression for cardinality of these sets:
\beq\label{Eq.Cardinality.K.n.r.s}
|\cK^s_{n,r}|=\bin{r}{s}\bin{n-1}{r-s-1}.
\eeq
If one thinks about $\cK_{n,r}$ as a simplex (see \S\ref{Sec.Knr.first.time}) then $\cK^s_{n,r}$ are collections of open cells/faces of codimension $s$.

 Some strata of $\cK_{n,r}$ are naturally isomorphic to (some strata of) $\cK_{n',r'}$ with lower $n',r'$. For example:
   $\cK^0_{n,r}\isom\cK_{n-r,r}$, by $(k_1,\dots,k_r)\to(k_1-1,\dots,k_r-1)$. Usually we identify $\cK_{n-r,r}$ with its image $\cK^0_{n,r}\sset\cK_{n,r}$.
    For example, we write: $|X\cap\cK^0_{n,r}|=|X\cap\cK_{n-r,r}|$.

    In the following we often split the sets into parts,  e.g. $\cK^s_{n,r}=\coprodl_{(i_1,\dots,i_s)\sset[1,\dots,r]}\cK^s_{n,r}(k_{i_1}=\cdots=k_{i_s}=0)$.
Here $\cK^s_{n,r}(k_{i_1}=\cdots=k_{i_s}=0)=\{\uk\in\cK^s_{n,r}\ |\ k_j=0$  iff  $j\in\{i_1,\dots,i_r\} \}$.
 As before, we get $\cK^s_{n,r}(k_{i_1}=\cdots=k_{i_s}=0)\isom\cK^0_{n,r-s}$.

 Sometimes we write this in the sloppy way $\cK^s_{n,r}=\coprodl_{(i_1,\dots,i_s)\sset[1,\dots,r]}\cK^0_{n,r-s}$. These formulas give immediate implications for the
 cardinality of the sets: $|\cK^s_{n,r}|=\bin{r}{s}|\cK^0_{n,r-s}|=\bin{r}{s}|\cK_{n-r+s,r-s}|$.

Note that starting from a subset $X\sseteq\cK_{n,r}$ we get the set $X\cap\cK^s_{n,r}$ and 
its subsets 
$X\cap \cK^s_{n,r}=\coprodl_{(i_1,\dots,i_s)\sset[1,\dots,r]}X\cap \cK^0_{n,r-s}=
\coprodl_{(i_1,\dots,i_s)\sset[1,\dots,r]}X\cap \cK_{n-r+s,r-s}$.
 Moreover, if the initial $X$ was $\Xi_r$-invariant and upward closed then so are all
the above subsets of type $X\cap K_{n-r+s,r-s}$.

\

{\bf Step 3.} Let a subset $X\sset\cK_{n,r}$ be  $\Xi_r$-invariant and upward closed, as above. Consider its averages over
 the strata, $\{Av_{\cK^s_{n,r}}X:=\frac{|X\cap\cK^s_{n,r}|}{|\cK^s_{n,r}|}\}_{s=0,\dots,r-1}$. We prove:
\beq
Av_{\cK^0_{n,r}}X\le Av_{\cK^1_{n,r}}X\le\cdots\le Av_{\cK^{r-1}_{n,r}}X
\eeq
and the equality holds iff $X=\cK_{n,r}$ or $X=\varnothing$.

The proof goes by  double induction on $(n,r)$. Recall that $\cK^s_{n,r}=\varnothing$ if $n<0$ or $r\le0$.
 Note that if $r=1$ or $n=1$ then the statement is empty ($\cK^s_{n,1}=\varnothing$ for $s>0$ and $\cK^s_{1,r}=\varnothing$ for $s<r-1$).
 For $r=2$ the proof is trivial.

Fix some $(n,r)$, suppose the statement holds for any $(n',r',X)$, with $n'<n$ and $r'<r$ and  $X$ an $\Xi_r$-invariant set which is upward closed.
 First we reduce all the inequalities to $Av_{\cK^0_{n,r}}X\le Av_{\cK^1_{n,r}}X$. Indeed:
\beq
Av_{\cK^s_{n,r}}X=\frac{|X\cap\cK^s_{n,r}|}{|\cK^s_{n,r}|}=
\suml_{(i_1,\dots,i_s)\sset[1,\dots,r]}\frac{|X\cap\cK^s_{n,r}(k_{i_1}=\cdots=k_{i_s}=0)|}{|\cK^s_{n,r}|}
\stackrel{\Xi_r\circlearrowright X}{=\joinrel=\joinrel=\joinrel=}
\bin{r}{s}\frac{|X\cap\cK^0_{n,r-s}|}{|\cK^s_{n,r}|}.
\eeq
(Here in the last equality we used the $\Xi_r$-invariance of $X$.)

In addition:
\beq
|X\cap\cK^1_{n,r+1-s}|=\suml^{r+1-s}_{j=1}|X\cap\cK^1_{n,r+1-s}(k_j=0)|\stackrel{\Xi_r\circlearrowright  X}{=\joinrel=\joinrel=\joinrel=}(r+1-s)|X\cap\cK^0_{n,r-s}|.
\eeq
Thus $Av_{\cK^s_{n,r}}X=\frac{\bin{r}{s}}{r+1-s}\frac{|X\cap\cK^1_{n,r+1-s}|}{|\cK^s_{n,r}|}$, while
 $Av_{\cK^{s-1}_{n,r}}X=\bin{r}{s-1}\frac{|X\cap\cK^0_{n,r+1-s}|}{|\cK^{s-1}_{n,r}|}$.
  Altogether:
\beq
Av_{\cK^s_{n,r}}X-Av_{\cK^{s-1}_{n,r}}X=\frac{\bin{r}{s}}{r+1-s}\frac{|\cK^1_{n,r+1-s}|}{|\cK^s_{n,r}|}\Big(
Av_{\cK^1_{n,r+1-s}}X-Av_{\cK^0_{n,r+1-s}}X\Big).
\eeq
(Here we used the equality of the coefficients:
$\frac{\bin{r}{s}}{r+1-s}\frac{|\cK^1_{n,r+1-s}|}{|\cK^s_{n,r}|}=
\bin{r}{s-1}\frac{|\cK^0_{n,r+1-s}|}{|\cK^{s-1}_{n,r}|}$.)

\

 If $s>1$ then  $Av_{\cK^1_{n,r+1-s}}X\ge Av_{\cK^0_{n,r+1-s}}X$ by the inductive assumption.
Thus it remains to prove that   $Av_{\cK^1_{n,r}}X\ge Av_{\cK^0_{n,r}}X$. Now we use the reduction:
\beqm
|X\cap\cK^1_{n,r}|=\suml^r_{j=1}|X\cap\cK^1_{n,r}(k_j=0)|=\suml^r_{j=1}|X\cap\cK^0_{n,r-1}|=
\suml^r_{j=1}|X\cap\cK_{n-r+1,r-1}|=\\=\suml^{r-2}_{s=0}\suml^r_{j=1}|X\cap\cK^s_{n-r+1,r-1}|
=\suml^{r-2}_{s=0}(s+1)|X\cap\cK^{s+1}_{n-r+1,r}|=\suml^{r-1}_{s=0}s |X\cap\cK^s_{n-r+1,r}|.
\end{multline}
Similarly:
\beq
|X\cap\cK^0_{n,r}|=\frac{1}{r}\suml^r_{j=1}|X\cap\cK_{n-r+1,r}(k_j>0)|=\frac{1}{r}\suml^r_{j=1}\suml^{r-1}_{s=0}|X\cap\cK^s_{n-r+1,r}(k_j>0)|=
\frac{1}{r}\suml^{r-1}_{s=0}(r-s)|X\cap\cK^s_{n-r+1,r}|.
\eeq
Therefore:
\beq
Av_{\cK^1_{n,r}}X-Av_{\cK^0_{n,r}}X=\frac{\suml^{r-1}_{s=0}s |X\cap\cK^s_{n-r+1,r}|}{|\cK^1_{n,r}|}-
\frac{\frac{1}{r}\suml^{r-1}_{s=0}(r-s)|X\cap\cK^s_{n-r+1,r}|}{|\cK^0_{n,r}|}=
\suml^{r-1}_{s=0}\underbrace{\Big(\frac{s}{|\cK^1_{n,r}|}-\frac{\frac{r-s}{r}}{|\cK^0_{n,r}|}\Big)|\cK^s_{n-r+1,r}|}_{\al_s}\underbrace{Av_{\cK^s_{n-r+1,r}}X}_{\be_s}.
\eeq
Here by the inductive assumption: $\be_{r-1}\ge\be_{r-2}\ge\cdots\ge\be_0\ge0$. By direct check:
\beq
\al_s=\frac{1}{\bin{n-2}{r-2}(n-1)}\bin{n-r}{r-s-1}\Big(n\bin{r-1}{s-1}-(r-1)\bin{r}{s}\Big).
\eeq
(We use the convention: $\bin{m}{n}=0$ if $n<0$ or $m<n$.)

We claim that $\suml^{r-1}_{s=k}\al_s=\frac{1}{\bin{n-2}{r-2}(n-1)}\bin{n-r}{r-k}\bin{r-1}{k}k$,  for $0\le k\le r-1$. (In particular, this sum is positive for $k>0$ and zero for $k=0$.)

 The case $k=0$ follows from
 $n\suml^{r-1}_{s=0}\bin{n-r}{r-s-1}\bin{r-1}{s-1}=(r-1)\suml^{r-1}_{s=0}\bin{n-r}{r-s-1}\bin{r}{s}$, which in turn follows from
the classical $\suml^p_{i=0}\bin{p}{i}\bin{q}{k-i}=\bin{p+q}{k}$.

Suppose the stated equality holds for some $k$, we prove it for $k+1$:
\beqm
\suml^{r-1}_{s=k+1}\al_s=\suml^{r-1}_{s=k}\al_s-\al_k=
\frac{\bin{n-r}{r-k}\bin{r-1}{k}k-\bin{n-r}{r-k-1}\Big(n\bin{r-1}{k-1}-(r-1)\bin{r}{k}\Big)}{\bin{n-2}{r-2}(n-1)}=\\
=\frac{\frac{(n-r)!}{(r-k)!(n-2r+k+1)!}\bin{r-1}{k}(k^2+k-2kr+r^2-r)}{\bin{n-2}{r-2}(n-1)}
=\frac{\bin{n-r}{r-k-1}\bin{r-1}{k+1}(k+1)}{\bin{n-2}{r-2}(n-1)},
\end{multline}
precisely as stated.

 Therefore:
\beqm
Av_{\cK^1_{n,r}}X-Av_{\cK^0_{n,r}}X=\suml^{r-1}_{s=0}\al_s\be_s=\underbrace{\al_{r-1}(\be_{r-1}-\be_{r-2})}_{\ge0}+
\underbrace{(\al_{r-1}+\al_{r-2})(\be_{r-2}-\be_{r-3})}_{\ge0}+
\\
+\cdots+\underbrace{ (\suml^{r-1}_{i=1}\al_i)(\be_1-\be_0)}_{\ge0}+\underbrace{(\suml^{r-1}_{i=0}\al_i)\be_0}_{=0}.
\end{multline}
As each term is non-negative we get: $Av_{\cK^1_{n,r}}X\ge Av_{\cK^0_{n,r}}X$.
Furthermore,  the equality holds iff $\{\be_{i+1}=\be_i\}_i$, i.e. either $\{X\cap\cK^i_{n,r}=\cK^i_{n,r}\}_i$ or $X=\varnothing$.

{\bf Step 4.} Finally, using $Av_{\cK^0_{n,r}}X\le Av_{\cK^1_{n,r}}X\le\cdots\le Av_{\cK^{r-1}_{n,r}}X$, we prove
$(Av_{\cK_{n,r}}X)(Av_{\cK_{n,r}}Y)\le Av_{\cK_{n,r}}(X\cap Y)$. As above, this is done by the double induction on $(n,r)$.
Note that the statement holds trivially for $r=1$ or $n=1$. Suppose it holds for any $(n',r')$ with $n'<n$ and $r'<r$.

First observe:
\beq
|\cK_{n,r}|\cdot Av_{\cK_{n,r}}(X\cap Y)=|X\cap Y\cap\cK_{n,r}|=\suml^{r-1}_{s=0}\suml_{(i_1,\dots,i_s)\sset[1,\dots,r]}
|X\cap Y\cap\cK^0_{n,r-s}(k_{i_1}=0=\cdots=k_{i_s})|.
\eeq
Further:
\beqm
|X\cap Y\cap\cK^0_{n,r-s}(k_{i_1}=0=\cdots=k_{i_s})|=|X\cap Y\cap\cK_{n-r+s,r-s}(k_{i_1}=0=\cdots=k_{i_s})|\ge
\\\\
\ge\frac{|X\cap\cK_{n-r+s,r-s}(k_{i_1}=0=\cdots=k_{i_s})|\cdot
|Y\cap\cK_{n-r+s,r-s}(k_{i_1}=0=\cdots=k_{i_s})|}{|\cK_{n-r+s,r-s}|}.
\end{multline}

The last inequality here is the induction assumption.

As was mentioned above, the cardinalities of the sets are related by: $|\cK_{n-r+s,r-s}|=|\cK^0_{n,r-s}|=\frac{|\cK^s_{n,r}|}{\bin{r}{s}}$.
 Similarly (as $X,Y$ are $\Xi_r$-invariant):
 $|X\cap\cK_{n-r+s,r-s}|=|X\cap\cK^0_{n,r-s}|=\frac{|X\cap\cK^s_{n,r}|}{\bin{r}{s}}$.
  Therefore:
\beq
|\cK_{n,r}|Av_{\cK_{n,r}}(X\cap Y)\ge\suml^{r-1}_{s=0}\suml_{(i_1,\dots,i_s)\sset[1,\dots,r]}
\frac{|X\cap\cK^s_{n,r}|\cdot|Y\cap\cK^s_{n,r}|}{\bin{r}{s}|\cK^s_{n,r}|}=
\suml^{r-1}_{s=0}\underbrace{\frac{|X\cap\cK^s_{n,r}|}{|\cK^s_{n,r}|}}_{\al_s}\cdot
\underbrace{\frac{|Y\cap\cK^s_{n,r}|}{|\cK^s_{n,r}|}}_{\be_s}\cdot
\underbrace{|\cK^s_{n,r}|}_{\ga_s}.
\eeq
Here we have $\al_0\le\al_1\le\cdots\le\al_{r-1}$ and $\be_0\le\be_1\le\cdots\le\be_{r-1}$, by {\em Step 3}, and $\ga_s>0$.

Thus we can use the following generalization of Chebyshev's sum inequality
\beq\label{eq:C<gen}
\big(\suml_s \ga_s\alpha_s\big)\big(\suml_{s'} \ga_{s'}\be_{s'}\big)\leq \big( \suml_{s'}\ga_{s'}\big)
\big(\sum _s\alpha_s\beta_s \ga_s\big),
\eeq
which basically is the summation $\suml_{s,s'}\ga_s\ga_{s'}(\al_s-\al_{s'})(\beta_s-\beta_{s'})\geq0$,
see \cite[p. 43]{Hardy-Littlewood-Polya}.

Hence:
\beqm
|\cK_{n,r}|^2\cdot Av_{\cK_{n,r}}(X\cap Y)\ge
\Big(\suml^{r-1}_{s'=0}\cK^{s'}_{n,r}\Big)\suml^{r-1}_{s=0}\underbrace{\frac{|X\cap\cK^s_{n,r}|}{|\cK^s_{n,r}|}}_{\al_s}\cdot
\underbrace{\frac{|Y\cap\cK^s_{n,r}|}{|\cK^s_{n,r}|}}_{\be_s}\cdot
\underbrace{|\cK^s_{n,r}|}_{\ga_s}\ge\\\ge
\Big(\suml^{r-1}_{s=0}|X\cap\cK^s_{n,r}|\Big)\Big(\suml^{r-1}_{s'=0}|Y\cap\cK^{s'}_{n,r}|\Big)
=|X\cap\cK_{n,r}|\cdot|Y\cap\cK_{n,r}|.
\end{multline}
This finishes the proof.
\epr

This theorem can be extended to the functions whose ``pushforwards'' are monotonic on $\quotients{\cK_{n,r}}{\Xi_r}$.
 Given $f:\cK_{n,r}\to\R_{\ge0}$ define $[f]:\quotients{\cK_{n,r}}{\Xi_r}\to\R_{\ge0}$ by: $[f]([x])=Av_{\Xi_r(x)}f$.
 (Note that this expression does not depend on the choice of the representative of $[x]$.)
\bcor
Given functions $f,g:\cK_{n,r}\to\R_{\ge0}$. Suppose $f$ is $\Xi_r$ invariant
 and $[f],[g]:\quotients{\cK_{n,r}}{\Xi_r}\to\R_{\ge0}$ are monotonic.
 Then $Av_{\cK_{n,r}}(f)\cdot Av_{\cK_{n,r}}(g)\le Av_{\cK_{n,r}}(fg)$, with equality iff one  of $[f]$, $[g]$ is constant.
\ecor
To prove this, define $\tg:\cK_{n,r}\to\R_{\ge0}$ by $\tg(x)=Av_{\Xi_r(x)}(g)$. Then
 $Av_{\cK_{n,r}}(g)=Av_{\cK_{n,r}}(\tg)$ and $Av_{\cK_{n,r}}(fg)=Av_{\cK_{n,r}}(f\tg)$.
 And by construction $f$, $\tg$ are monotonic on $\cK_{n,r}$.
  Now, apply the theorem.

 \

 \

 We can reformulate theorem \ref{Thm.Averaging.over.Knr} to compare the geometric averages.
\bcor\label{Thm.Inequality.Geometric.Averages}
Given two non-decreasing functions, $f:\cK_{n,r}\to\R_{\ge1}$ and $g:\cK_{n,r}\to\R_{\ge0}$.
Suppose at least one of them is $\Xi_r$-invariant.
Then \[\Big(Av^G_{\cK_{n,r}}(f)\Big)^{Av_{\cK_{n,r}}(g)}\ge Av^G_{\cK_{n,r}}(f^g).\]
 \ecor
To prove this take the logarithm of both sides. Note that $ln(Av^G_{\cK_{n,r}}(f))=Av_{\cK_{n,r}}(ln(f))$
where  $ln(f):\cK_{n,r}\to\R_{\ge0}$ is still a  non-decreasing function.
Similarly $ln(Av^G_{\cK_{n,r}}(f^g))=Av_{\cK_{n,r}}(g\cdot ln(f))$.
Now, apply the theorem.

\section{An application to mixed (co-)volumes}

\subsection{Newton polyhedra}\label{Sec.Background.Newton.Diagrams}
Let $f(x_1,\dots,x_N)=\suml_I a_I \ux^I$ be a  power series over some field $\k$. Consider the support
 of its monomials, $Supp(f):=\{I\in\Z^N_{\ge0}|\ a_I\neq0\}$. The Newton polyhedron is defined as
 the convex hull,
\beq
\Ga^+=\Ga^+_f:=Conv(Supp(f)+\R^N_{\ge0}).
\eeq
 It has compact faces and unbounded faces.
The Newton polyhedron is called convenient if $\Ga^+$ intersects all
 the coordinate axes (i.e. $f$ contains monomials of type $x^{m_1}_1,\dots,x^{m_N}_N$).
We always assume  $\Ga^+$ to be convenient.
(Though in the present note we do not study the analytic properties of the power series
related with their Newton diagrams, by the above definition we wish to emphasize the
main motivation supported by algebraic geometry and complex analysis.)

\subsection{Mixed (co)volumes}\label{Sec.Background.Mixed.Covolumes}
Given several convex bodies $A_1,\dots,A_r$ in $\R^N$, consider their scaled Minkowski sum, $\la_1 A_1+\cdots+\la_r A_r$.
 The mixed volumes are defined as the coefficients in the expansion:
\beq
Vol_N(\la_1 A_1+\cdots+\la_r A_r)=\suml_{\uk\in\cK_{N,r}}\bin{N}{k_1,\dots,k_r}Vol\Big((A_1)^{k_1},\dots,(A_r)^{k_r}\Big)(\prod^r_{i=1}\la^{k_i}_i).
\eeq
Dually, given a convenient Newton polyhedron, $\Ga^+\sset\R^N_{\ge0}$, consider its covolume, i.e. the volume of the complement: $coVol(\Ga^+):=Vol_N(\R^N_{\ge0}\smin\Ga^+)$.
Given a collection of Newton polyhedra, $\{\Ga^+_i\}_{i=1}^r$,
consider their scaled Minkowski sum,
$\la_1\Ga^+_1+\cdots+\la_r\Ga^+_r$.
The covolume of this sum is a polynomial in $\{\la_i\}$, see e.g. \cite[Theorem 10.4]{Kaveh-Khovanskii-Mix.Mult}:
\beq
coVol(\la_1\Ga^+_1+\cdots+\la_r\Ga^+_r)=
\suml_{\uk\in\cK_{N,r}}\bin{N}{k_1,\dots,k_r}coVol\Big((\Ga^+_1)^{k_1},\dots,(\Ga^+_r)^{k_r}\Big)(\prod^r_{i=1}\la^{k_i}_i).
\eeq
The mixed covolumes are the (positive) coefficients $coVol\Big((\Ga^+_1)^{k_1},\dots,(\Ga^+_r)^{k_r}\Big)$.

Here $coVol\Big((\Ga^+_1)^{k_1},\dots,(\Ga^+_r)^{k_r}\Big)$ is a shorthand for
$coVol\Big(\underbrace{\Ga^+_1,\dots,\Ga^+_1}_{k_1},\dots,\underbrace{\Ga^+_r,\dots,\Ga^+_r}_{k_r}\Big)$, for $k_1+\cdots+k_r=N$,
 similarly to $Vol\Big((A_1)^{k_1},\dots,(A_r)^{k_r}\Big)$.

We list some of the basic properties of the mixed (co)volumes, see e.g. \cite[\S2]{Kaveh-Khovanskii-Conv.Bod}:
\li They are symmetric and multilinear, e.g. $coVol(\Ga^+_{11}+\Ga^+_{12},\Ga_2^+,\dots,\Ga^+_N)=coVol(\Ga^+_{11},\Ga_2^+,\dots,\Ga^+_N)+
coVol(\Ga^+_{12},\Ga_2^+,\dots,\Ga^+_N)$.
\li If $\Ga^+_i=d_i\Ga^+$ for any $i=1,\ldots ,r$ then
\beq
coVol_N(\suml_i\la_i\Ga^+_i)=coVol_N(\suml_i \la_id_i\Ga^+)=(\suml_i \la_id_i)^N coVol_N(\Ga^+)=
\suml_{\uk\in\cK_{n,r}}\bin{N}{k_1,\dots, k_r}(\prodl^r_{i=1} (\la_i d_i)^{k_i}) coVol_N(\Ga^+).
\eeq
A similar statement holds  for the convex bodies: $Vol_N(\suml_i\la_iA_i)=Vol_N(\suml_i \la_id_iA)=(\suml_i d_i\la_i)^N Vol_N(A)$.
\li The mixed volumes of the convex bodies in $\R^N$ satisfy the Alexandrov-Fenchel inequality, \cite{Alexandrov-1937}, \cite{Alexandrov-1938}, \cite{Fenchel}:
\beq\label{Eq.Alex.Fench.ineq}
Vol(A_1,\dots,A_N)^2\ge Vol(A_1,A_1,A_3\dots,A_N)Vol(A_2,A_2,A_3\dots,A_N).
\eeq
The dual property of the mixed covolumes of Newton polyhedra  was proved in \cite{Teissier1978}, \cite{Rees-Sharp-1978}, \cite{Katz1988}, \cite[Appendix]{Teissier2004},
\cite[Theorem 10.5]{Kaveh-Khovanskii-Mix.Mult}:
\beq\label{Eq.Teissier.Ineq} coVol(\Ga^+_1,\Ga^+_2,\cdots,\Ga^+_N)^2\le
coVol(\Ga^+_1,\Ga^+_1,\Ga^+_3,\cdots,\Ga^+_N)coVol(\Ga^+_2,\Ga^+_2,\Ga^+_3,\cdots,\Ga^+_N).
\eeq

\subsection{A generalization of the Alexandrov-Fenchel-Teissier inequalities}\label{Sec.Alexandrov.Fenchel.Teissier.ineqs}
The inequalities of equations \eqref{Eq.Alex.Fench.ineq} and \eqref{Eq.Teissier.Ineq} correspond to the triples of points on a segment in $\cK_{n,r}$: $\{k_1+2,k_2,k_3,\dots,\}$, $\{k_1+1,k_2+1,k_3,\dots\}$,
$\{k_1,k_2+2,k_3,\dots\}$.
One would ask for the general property, the Jensen-type inequality:

{\em Suppose a collection of lattice points $\uk^{(1)},\uk^{(2)},\dots,\uk^{(s)}\in\cK_{N,r}$, satisfy:
 $\uk:=\frac{\uk^{(1)}+\cdots+\uk^{(s)}}{s}\in\cK_{N,r}$ (i.e. $\uk$ defined in this way is a lattice point too).
 Then
 $coVol\big((\Ga^+)^{\uk}\big)^s\le \prodl^s_{i=1} coVol\big((\Ga^+)^{\uk^{(i)}}\big)$ and $Vol\big(\uA^{\uk}\big)^s\ge \prodl^s_{i=1} Vol\big(\uA^{\uk^{(i)}}\big)$.}

In particular, Gromov in \cite{gromov1990convex} asked whether this holds at least for  points
$\uk^{(1)},\uk^{(2)},\dots,\uk^{(s)}$ sitting inside a low dimensional linear subspace of $\cK_{n,r}$.
As it was shown by Burda in \cite{Burda} such a property in general fails, e.g. he provided
 an example with
\beq\label{Eq.Burda}
Vol_3(A_1,A_2,A_3)^3< Vol_3(A_1,A_1,A_2)Vol_3(A_2,A_2,A_3)Vol_3(A_3,A_3,A_1).
\eeq
Yet, Theorem \ref{Thm.Averaging.over.Knr} and Corollary \ref{Thm.Inequality.Geometric.Averages} allow
us to establish  the  `corrected version' of the above question
(see Corollary \ref{Thm.Averaging.Mixed.Volumes}).

\

Given a pair $(n',r')$ fix some $1\le r\le r'$ and denote $n:=n'+r-r'$.
Fix some Newton polyhedra, $\Ga^+_1,\dots,\Ga^+_r,\Ga^+_{r+1},\dots,\Ga^+_{r'}\sset\R^{n'}$.
 Define the function $Mix.coVol:$ $\cK_{n,r}\to\R_{\ge0}$
 by $\uk\to coVol\Big((\Ga^+_1)^{k_1},\dots,(\Ga^+_r)^{k_r},\Ga^+_{r+1},\dots,\Ga^+_{r'}\Big)$.
   Fix some convex bodies $A_1,\dots,A_r,A_{r+1},\dots,A_{r'}\sset\R^{n'}$ and define the function
\beq
Mix.Vol:\ \cK_{n,r}\to\R_{\ge0},\quad \uk\to Vol\Big(A_1^{k_1},\dots,A_r^{k_r},A_{r+1},\dots,A_{r'}\Big).
\eeq
 These functions correspond to the particular embedding $\cK_{n,r}\hookrightarrow\cK_{n',r'}$
 by $(k_1,\dots,k_r)\to(k_1,\dots,k_r,1,\dots,1)$.
\bcor\label{Thm.Averaging.Mixed.Volumes}
Let $C:\cK_{n,r}\to\R_{\ge0}$ be a non-decreasing $\Xi_r$-invariant function.
\\1.  $Av_{\cK_{n,r}}(C)\cdot Av_{\cK_{n,r}}(Mix.coVol)\le Av_{\cK_{n,r}}(C\cdot Mix.coVol)$.
\\2.  Suppose all the values of $Mix.coVol$ and $Mix.Vol$ are $\ge1$. Then
$\Big(Av^G_{\cK_{n,r}}(Mix.coVol)\Big)^{Av_{\cK_{n,r}}(C)}\le Av^G_{\cK_{n,r}}(Mix.coVol^C)$ and
$\Big(Av^G_{\cK_{n,r}}(Mix.Vol)\Big)^{Av_{\cK_{n,r}}(C)}\ge Av^G_{\cK_{n,r}}(Mix.Vol^C)$.
\ecor
\bpr {\bf Part 1.}

{\em Step 1. } First, using the inequality $\sqrt{ab}\le\frac{a+b}{2}$ one gets from equation \eqref{Eq.Teissier.Ineq} the weaker convexity property:
\beq\label{Eq.Mixed.Vol.convexity.weak}
coVol(\Ga^+_1,\Ga^+_2,\cdots,\Ga^+_N)\le
\frac{coVol(\Ga^+_1,\Ga^+_1,\Ga^+_3,\cdots,\Ga^+_N)+coVol(\Ga^+_2,\Ga^+_2,\Ga^+_3,\cdots,\Ga^+_N)}{2}.
\eeq
The function $Mic.coVol$ is not $\Xi_r$-invariant, for $\Xi_r$ acting on the first $r$ indices. Consider its symmetrization:
\beq
 [Mix.coVol]:\ \cK_{n,r}\to\R_{\ge0}, \quad\quad
 \uk\to Av_{\Xi_r(\uk)}(Mix.coVol)=\frac{\suml_{\si\in\Xi_r}coVol\Big((\Ga^+_{\si(1)})^{k_1},\dots,(\Ga^+_{\si(r)})^{k_r},\Ga^+_{r+1},\dots,\Ga^+_{r'}\Big)}{|\Xi_r(\uk)|}.
\eeq
We prove that this symmetrization is a non-decreasing function, i.e.
\label{Thm.Mixed.Volume.Non.Decreasing.Func}
if $\uk\succeq\uk'\in\cK_{n,r}$ then $[Mix.coVol](\uk)\ge [Mix.coVol](\uk')$.

It is enough to check the inequality for elementary steps: if $k_1+1\le k_2-1$ then
$[Mix.coVol]([k_1,k_2,k_3,\dots])\ge [Mix.coVol]([k_1+1,k_2-1,k_3,\dots])$.

Suppose $k_1+k_2=2l\in2\Z$ . Then the convexity property of the mixed volumes reads:
\beqm
[Mix.coVol](l,l,k_3,\dots)=
\suml_{\si\in\Xi_r}\frac{coVol\Big((\Ga^+_{\si(1)})^{l},(\Ga^+_{\si(2)})^{l},(\Ga^+_{\si(3)})^{k_3}\dots,(\Ga^+_{\si(r)})^{k_r},\dots\Big)}{r!}
\stackrel{eq. \eqref{Eq.Mixed.Vol.convexity.weak}}{\le}
\\\le
\suml_{\si\in\Xi_r}\frac{coVol\Big((\Ga^+_{\si(1)})^{l-1},(\Ga^+_{\si(2)})^{l+1},(\Ga^+_{\si(3)})^{k_3}\dots,(\Ga^+_{\si(r)})^{k_r},\dots\Big)+
coVol\Big((\Ga^+_{\si(1)})^{l+1},(\Ga^+_{\si(2)})^{l-1},(\Ga^+_{\si(3)})^{k_3}\dots,(\Ga^+_{\si(r)})^{k_r},\dots\Big)}{2r!}=
\\=
\suml_{\si\in\Xi_r}\frac{coVol\Big((\Ga^+_{\si(1)})^{l-1},(\Ga^+_{\si(2)})^{l+1},(\Ga^+_{\si(3)})^{k_3}\dots,(\Ga^+_{\si(r)})^{k_r},\dots\Big)}{r!}
=[Mix.coVol](l-1,l+1,k_3,\dots).
\end{multline}
Similarly,  $[Mix.coVol](l-1,l+1,k_3,\dots)\le\frac{[Mix.coVol](l,l,k_3,\dots)+[Mix.coVol](l-2,l+2,k_3,\dots)}{2}$.
 Combining with the previous one it gives $[Mix.coVol](l-2,l+2,k_3,\dots)\ge [Mix.coVol](l-1,l+1,k_3,\dots)$.
  One continues in the same way.

  The case $k_1+k_2=2l+1\in2\Z+1$ is treated similarly.

 {\em Step 2.} Now we have two $\Xi_r$-invariant, non-decreasing functions, then theorem \ref{Thm.Averaging.over.Knr} gives:
\beq
Av_{\cK_{n,r}}(C)Av_{\cK_{n,r}}([Mix.coVol])\le Av_{\cK_{n,r}}(C\cdot[Mix.coVol]).
\eeq
 It remains to add that
  $Av_{\cK_{n,r}}([Mix.coVol])=Av_{\cK_{n,r}}(Mix.coVol)$ and $Av_{\cK_{n,r}}(C\cdot[Mix.coVol])=Av_{\cK_{n,r}}(C\cdot Mix.coVol)$.

\

{\bf Part 2.}   Define the function $f:\ \cK_{n,r}\to\R_{\ge0}$  by $f(\uk)=ln\Big(coVol\big(A_1^{k_1},\dots,A_r^{k_r},A_{r+1},\dots,A_{r'}\big)\Big)$. Equation \eqref{Eq.Teissier.Ineq} implies:
$f(k_1,k_2,k_3,\dots)\le \frac{f(k_1,k_1,k_3,\dots)+f(k_2,k_2,k_3,\dots)}{2}$. Then, repeat Part 1 for the functions $f$, $C$ to get:
\beq
Av_{\cK_{n,r}}\Big(ln\big(coVol(\dots)^{Av_{\cK_{n,r}}(C)}\big)\Big)=Av_{\cK_{n,r}}(C)Av_{\cK_{n,r}}(f)\le Av_{\cK_{n,r}}(C\cdot f)=Av_{\cK_{n,r}}\Big(ln\big(coVol(\dots)^{C(\dots)}\big)\Big).
\eeq
Now take the exponent.

The proof of the inequality for mixed volumes is the same, it is based now on equation \eqref{Eq.Alex.Fench.ineq}.
  \epr

\bex
$\bullet$ The case $r=1$ is `empty' as $\cK_{n,1}=\{n\}$.
\li As the simplest case, suppose $r=2\le r'$, $n=2\le n'$. Then part 2 of the corollary states just the ordinary Alexandrov-Fenchel/Teissier inequalities.
\li An extremal case is $r=r'$, then $n=n'$. For $n=r=3$ fix the numbers $C_{(3,0,0)}=a\ge C_{(2,1,0)}=b\ge C_{(1,1,1)}=c\ge0$. Extend this to the function $C:\cK_{3,3}\to\R_{\ge0}$ by $\Xi_3$-action.
 We get a $\Xi_3$-invariant non-decreasing function. Then Part 2 of the proposition gives:
\beq
Vol(A_1,A_2,A_3)^{3a+6b-9c}\ge \Big(\prodl^3_{i=1}Vol(A_i,A_i,A_i)\Big)^{7a-6b-c}
\Big(\prodl_{j\neq i}Vol(A_i,A_i,A_j)\Big)^{4b-3a-c}
\eeq
(This formula becomes equality in the simplest case $\{A_i=\la_i A\}$. Then $Vol_3(A_i,A_i,A_i)=d^3_i Vol_3(A)$, $Vol_3(A_i,A_i,A_j)=d^2_id_j Vol_3(A)$, $Vol_3(A_1,A_2,A_3)=d_1 d_2d_3 Vol_3(A)$.)
\eex
\bex
One can get a bigger class of inequalities by multiplying the inequalities of part 2 of the corollary for different embeddings $\cK_{n,r}\hookrightarrow\cK_{n',r'}$. (We do not know whether they contain/generate all the possible inequalities on mixed volumes.)
For example, fix $r=2=n<r'=3=n'$. Consider the embeddings
 $(k_1,k_2)\to(k_1,k_2,1)$, $(k_1,k_2)\to(1,k_1,k_2)$, $(k_1,k_2)\to(k_1,1,k_2)$. Each of them gives the ordinary Aleandrov-Fenchel inequality. The product of these inequalities is:
\beqm
Vol_3(A_1,A_2,A_3)^6\ge\\ Vol_3(A_1,A_1,A_2)Vol_3(A_1,A_2,A_2)Vol_3(A_2,A_2,A_3)Vol_3(A_2,A_3,A_3)Vol_3(A_3,A_3,A_1)
Vol_3(A_3,A_1,A_1).
\end{multline}
This can be considered as a `corrected' version of equation \eqref{Eq.Burda}.
\eex

\bex
In \cite{Kerner-Nemethi.Durfee.Nnd} we use this corollary in the particular case of $C(\uk)=\bin{n+r}{k_1+1,\dots,k_r+1}:=\frac{(n+r)!}{(k_1+1)!\cdots(k_r+1)!}$.
 (Note that this expression is $\Xi_r$-invariant and non-decreasing on $\cK_{n,r}$.)
We claim that:
\beqm\label{Eq.Ineq.To.Prove}
\Big(\suml_{\uk\in\cK_{n,r}}\bin{n+r}{k_1+1,\dots,k_r+1}\Big)
\suml_{\substack{\uk\in\cK_{n+r,r}\\k_1,\dots,k_r\ge1}}coVol_{n+r}\Big((\Ga^+_1)^{k_1},\dots,(\Ga^+_r)^{k_r}\Big)
\ge\\
\ge
\bin{n+r-1}{n}\suml_{\substack{\uk\in\cK_{n+r,r}\\k_1,\dots,k_r\ge1}}
\bin{n+r}{k_1+1,\dots,k_r+1}coVol_{n+r}\Big((\Ga^+_1)^{k_1},\dots,(\Ga^+_r)^{k_r}\Big)
\end{multline}
and equality occurs iff $\Ga^+_1=\cdots=\Ga^+_r$.

Note that $\cK_{n+r,r}\cap\{k_1,\dots,k_r\ge1\}=\cK^0_{n+r,r}\isom\cK_{n,r}$ in the notations of the proof of theorem \ref{Thm.Averaging.over.Knr}. Therefore the inequality is the
 comparison of averages of the numbers
  $\{coVol_{n+r}\Big((\Ga^+_1)^{k_1+1},\dots,(\Ga^+_r)^{k_r+1}\Big)\}_{\uk\in\cK_{n,r}}$ and
$\{\bin{n+r}{k_1+1,\dots,k_r+1}\}_{\uk\in\cK_{n,r}}$.

Note that $\{\bin{n+r}{k_1+1,\dots,k_r+1}$ is symmetric under permutations of $\uk$ and gives
 a non-increasing function on $\quotients{\cK_{n,r}}{\Xi_r}$.

Thus, corollary \ref{Thm.Averaging.Mixed.Volumes} gives:
\beqm
Av_{\cK_{n,r}}\Big\{\bin{n+r}{k_1+1,\dots,k_r+1}\Big\}
Av_{\cK_{n,r}}\Big\{coVol_{n+r}\Big((\Ga^+_1)^{k_1},\dots,(\Ga^+_r)^{k_r}\Big)\Big\}\ge \\\ge
Av_{\cK_{n,r}}\Big\{\bin{n+r}{k_1+1,\dots,k_r+1}coVol_{n+r}\Big((\Ga^+_1)^{k_1},\dots,(\Ga^+_r)^{k_r}\Big)\Big\}
\end{multline}
which is precisely the inequality \eqref{Eq.Ineq.To.Prove}.

Moreover, as $\bin{n+r}{k_1+1,\dots,k_r+1}$ is not constant on $\quotients{\cK_{n,r}}{\Xi_r}$, the equality
 occurs iff the symmetrization of the function $coVol_{n+r}\Big((\Ga^+_1)^{k_1},\dots,(\Ga^+_r)^{k_r}\Big)$ is constant on $\cK_{n,r}$; which means $\Ga^+_1=\cdots=\Ga^+_r$.
\eex

\end{document}